\documentclass[11pt]{article}

\usepackage{graphicx,psfrag}
\usepackage{psfrag,eepic,epsfig}
\usepackage{sectsty}
\usepackage{setspace}
\usepackage{color}

\usepackage{setspace}
\usepackage{amsmath, epsfig, cite, ifpdf,float}
\usepackage{amssymb}
\usepackage{amsfonts}
\usepackage{latexsym}
\usepackage{graphicx,psfrag,eepic,rawfonts}
\usepackage{enumerate}
\usepackage{indentfirst}
\input{epsf}

\newtheorem{thm}{Theorem}[section]

\newtheorem{lem}[thm]{Lemma}

\numberwithin{equation}{section}

\makeatletter \@addtoreset{equation}{section} \makeatother

\begin{document}
\rule{0cm}{1cm}

\begin{center}
{\Large\bf Tricyclic graphs with maximal revised Szeged index
}
\end{center}

\begin{center}
{\small Lily Chen, Xueliang Li, Mengmeng Liu\\
Center for Combinatorics, LPMC\\
Nankai University, Tianjin 300071, China\\
Email:  lily60612@126.com, lxl@nankai.edu.cn, liumm05@163.com}
\end{center}

\begin{center}
\begin{minipage}{120mm}
\begin{center}
{\bf Abstract}
\end{center}

{\small The revised Szeged index of a graph $G$ is defined as
$Sz^*(G)=\sum_{e=uv \in E}(n_u(e)+ n_0(e)/2)(n_v(e)+ n_0(e)/2),$
where $n_u(e)$ and $n_v(e)$ are, respectively, the number of
vertices of $G$ lying closer to vertex $u$ than to vertex $v$ and
the number of vertices of $G$ lying closer to vertex $v$ than to
vertex $u$, and $n_0(e)$ is the number of vertices equidistant to
$u$ and $v$. In this paper, we give an upper bound of the revised
Szeged index for a connected tricyclic graph, and also characterize
those graphs that achieve the upper bound.
}

\vskip 3mm

\noindent {\bf Keywords:} Wiener index, Szeged index, Revised Szeged
index, tricyclic graph.

\vskip 3mm


\end{minipage}
\end{center}

\section{Introduction}

All graphs considered in this paper are finite, undirected and
simple. We refer the reader to \cite{bm} for terminology and
notation not given here. Let $G$ be a connected graph with vertex
set $V(G)$ and edge set $E(G)$. For $u,v \in V(G)$, $d_G(u,v)$
denotes the {\it distance} between $u$ and $v$ in $G$,  we use
$d(u,v)$ for short, if there is no ambiguity. The {\it Wiener index}
of $G$ is defined as
$$
W(G)=\displaystyle\sum_{\{u,v\}\subseteq V(G)} d_G(u,v).
$$
This topological index has been extensively studied in the
mathematical literature; see, e.g., \cite{GSM,GYLL}. Let $e=uv$ be
an edge of $G$, and define three sets as follows:
$$
N_u(e) = \{w \in V(G): d_G(u,w)< d_G(v,w)\},
$$
$$
N_v(e) = \{w \in V(G): d_G(v,w)< d_G(u,w)\},
$$
$$
N_0(e) = \{w \in V(G): d_G(u,w)=d_G(v,w)\}.
$$
Thus, $\{N_u(e),N_v(e),N_0(e)\}$ is a partition of the vertices of
$G$ respect to $e$. The number of vertices of $N_u(e)$, $N_v(e)$ and
$N_0(e)$ are denoted by $n_u(e)$, $n_v(e)$ and $n_0(e)$,
respectively. A long time known property of the Wiener index is the
formula \cite{GP,W}:
$$
W(G) = \displaystyle\sum_{e=uv \in E(G)} n_u(e) n_v(e),
$$
which is applicable for trees. Motivated by the above formula,
Gutman \cite{G} introduced a graph invariant, named as the {\it
Szeged index}, as an extension of the Wiener index and defined by
$$
Sz(G) = \displaystyle\sum_{e=uv \in E(G)}n_u(e) n_v(e).
$$
Randi\'c \cite{R} observed that the Szeged index does not take into
account the contributions of the vertices at equal distances from
the endpoints of an edge, and so he conceived a modified version of
the Szeged index which is named as the {\it revised Szeged index}.
The revised Szeged index of a connected graph $G$ is defined as
$$
Sz^*(G) = \displaystyle\sum_{e=uv \in E(G)}\left(n_u(e)+
\frac{n_0(e)}{2}\right)\left(n_v(e)+ \frac{n_0(e)}{2}\right).
$$

Some properties and applications of these two topological indices
have been reported in \cite{CLL,arX,auto,I,LL,PR,PZ,SGB}. In \cite{AH},
Aouchiche and Hansen showed that for a connected graph $G$ of order
$n$ and size $m$, an upper bound of the revised Szeged index of $G$
is $\frac{n^2m}{4}$. In \cite{XZ}, Xing and Zhou determined the
unicyclic graphs of order $n$ with the smallest and the largest
revised Szeged indices for $n\geq 5$, and they also determined the
unicyclic graphs of order $n$ with the unique cycle of length $r \
(3\leq r\leq n)$, with the smallest and the largest revised Szeged
indices. In \cite{LL}, we identified those graphs whose revised
Szeged index is maximal among bicyclic graphs. In this paper, we
give an upper bound of the revised Szeged index for a connected
tricyclic graph, and also characterize those graphs that achieve the
upper bound.

\begin{figure}[h,t,b,p]
\begin{center}
\scalebox{0.9}[0.9]{\includegraphics{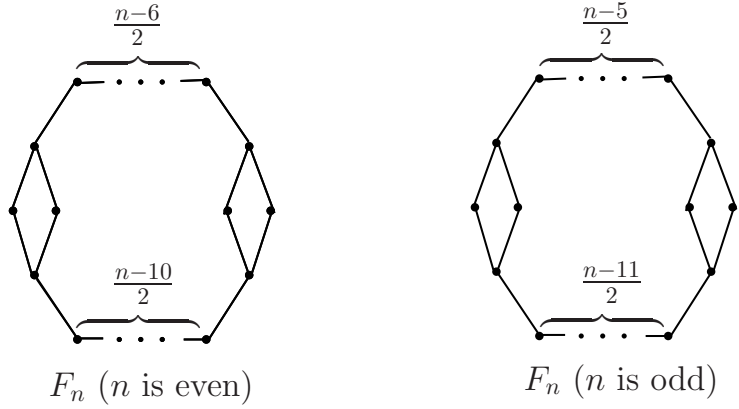}}
\end{center}
\caption{The graph for Theorem \ref{th1}}\label{f1}
\end{figure}

\begin{thm} \label{th1}
Let $G$ be a connected tricyclic graph
$G$ of order $n \geq 29$. Then
$$
Sz^*(G)\leq \left\{
\begin{array}{ll}
(n^3+2n^2-16)/4,& \mbox {if $n$ is even},\\
(n^3+2n^2-18)/4, & \mbox {if $n$ is odd}.
\end{array}
\right.
$$
with equality if and only if $G \cong F_n$ (see Figure \ref{f1}).
\end{thm}

\section{Main result}

It is easy to check that
$$
Sz^*(F_n) =  \left\{
\begin{array}{ll}
(n^3+2n^2-16)/4,& \mbox {if $n$ is even},\\
(n^3+2n^2-18)/4, & \mbox {if $n$ is odd}.
\end{array}
\right.
$$
i.e., $F_n$ satisfies the equality of Theorem \ref{th1}.

So, we are left to show that for any connected tricyclic graph $G_n$
of order $n \geq 29$, other than $F_n$, $Sz^*(G_n)<Sz^*(F_n)$. Using the
fact that $n_u(e)+ n_v(e)+n_0(e)=n$ and $m=n+2$, we have
\begin{eqnarray*}
Sz^*(G) & = & \displaystyle \sum_{e=uv \in E(G)} \left(n_u(e)+
\frac{n_0(e)}{2}\right)\left(n_v(e)+ \frac{n_0(e)}{2}\right) \\
&=& \displaystyle \sum_{e=uv \in E(G)}
\left(\frac{n+n_u(e)-n_v(e)}{2}\right)
\left(\frac{n-n_u(e)+n_v(e)}{2}\right)\\
 &= & \displaystyle \sum_{e=uv \in E(G)}
 \frac{n^2-(n_u(e)-n_v(e))^2}{4} \nonumber\\
 &= & \frac{mn^2}{4}-\frac{1}{4}\displaystyle \sum_{e=uv \in
 E(G)}(n_u(e)-n_v(e))^2.\\
 &=& \frac{n^3+2n^2}{4}-\frac{1}{4}\displaystyle \sum_{e=uv \in
 E(G)}(n_u(e)-n_v(e))^2
\end{eqnarray*}
For convenience, let $\delta (e)= |n_u(e)-n_v(e)|$, where $e=uv.$ We have
$$ Sz^*(G)  =
\frac{n^3+2n^2}{4}-\frac{1}{4}\displaystyle \sum_{e=uv \in
 E(G)} \delta ^2(e)   \eqno(1)
$$

\subsection{Proof for tricyclic graphs with connectivity 1}

\begin{lem}\label{lem1}
Let $G$ be a connected tricyclic graph of order $n\geq 12$ with at
least one pendant edge. Then
$$
Sz^*(G_n)<Sz^*(F_n)
$$
\end{lem}

\begin{pf}
Let $e'=xy$ be a pendant edge and $d(y)=1$. Then, for $n \geq 12,$ we
have
\begin{eqnarray*}
\displaystyle\sum_{e=uv\in E}(n_u(e)-n_v(e))^2 & \geq &
(n_x(e')-n_y(e'))^2\\
&=& (n-1-1)^2 \\
&>& 18.
\end{eqnarray*}
Combining with equality $(1)$, the result follows.
\end{pf}
\begin{qed}
\end{qed}

\begin{lem}\label{lem2.1}
Let $G$ be a connected tricyclic graph of order $n\geq 12$ without
pendant edges but with a cut vertex. Then, we have
$$
Sz^*(G)<Sz^*(F_n)
$$
\end{lem}

\begin{pf}
Suppose that $u$ is a cut vertex. Since $G$ is a tricyclic graph
without pendant edge, $G$ is composed of a bicyclic graph
 $B$ and a cycle $C$ and $V(B)\cap V(C)=\{u\}$.
 It is obvious that $|V(B)| \geq 4$. If $C$ is even, for
every edge $e$ in $C$, we have $\delta(e)=|V(B)|-1=n-|V(C)|.$ So
\begin{eqnarray*}
\displaystyle\sum_{e\in E(G)}\delta^2(e)  \geq
\displaystyle\sum_{e\in E(C)}\delta^2(e) =
 |E(C)|(|V(B)|-1)^2 \geq 4 \times 3^2 >18.
\end{eqnarray*}
If $C$ is odd, for all edges in $C$ but the edge $xy$ such
that $d(u,x)=d(u,y)$, we have
$\delta(e)=|V(B)|-1=n-|V(C)|.$ So
\begin{eqnarray*}
\displaystyle\sum_{e\in E(G)}\delta^2(e)  \geq
\displaystyle\sum_{e\in E(C)}\delta^2(e) =
 (|E(C)|-1)(|V(B)|-1)^2 .
\end{eqnarray*}
If $|E(C)| \geq 5$, then $\displaystyle\sum_{e\in E(G)}\delta^2(e)  > 18.$
If $|E(C)| = 3$, then $|V(B)|-1=n-|V(C)| \geq 9$, so
$\displaystyle\sum_{e\in E(G)}\delta^2(e) > 18$.

Combining with equality $(1)$, this completes the proof.
\end{pf}
\begin{qed}
\end{qed}

\subsection{Proof for 2-connected tricyclic graphs}
\begin{figure}[h,t,b,p]
\begin{center}
\scalebox{0.8}[0.8]{\includegraphics{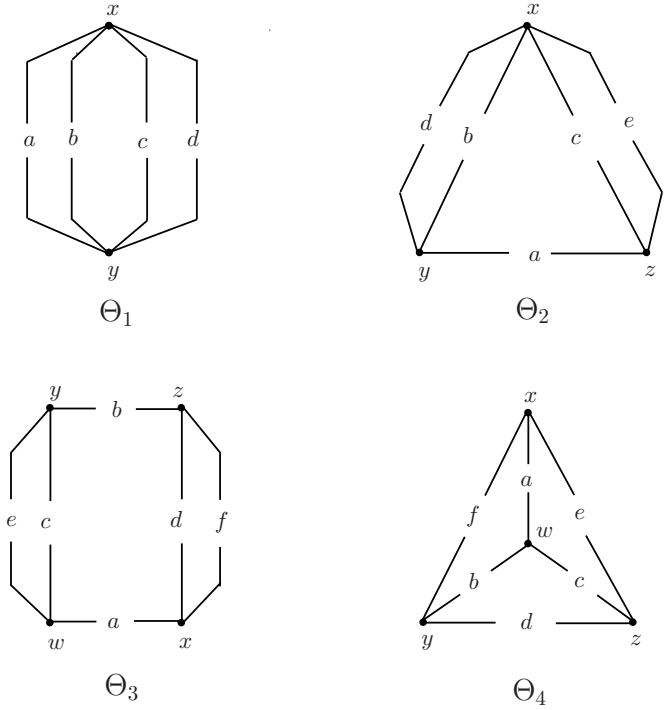}}
\end{center}
\caption{Four cases for 2-connected tricyclic graphs.}\label{f2}
\end{figure}

In this section, $\kappa(G)\geq 2$, then it must be one of the graphs
depicted in Figure \ref{f2}. The letters $a,b,\ldots,f$ stand  for the lengths
of the corresponding paths between vertices of degree greater than 2.
For the sake of brevity, we refer to these paths as $P(a),P(b),\ldots,
P(f)$, respectively. In the statement of the following lemmas, we call
these four graphs in Figure \ref{f2}  as $\Theta_1, \Theta_2,\Theta_3$ and $\Theta_4$,
respectively.

\begin{lem}\label{lem3.1}
Let $G$ be a $\Theta_1$-graph composed of four paths $P_1$,
$P_2,P_3$ and $P_4$, and $e=uv \in E(G)$. Then $|n_u(e)-n_v(e)| \leq 1$ if and
only if $e$ is in the middle of an odd path of the four paths
$P_1$, $P_2$, $P_3$ and $P_4$.
\end{lem}

\begin{pf}
Assume that
$e=uv$ belongs to $P_i \ (1\leq i\leq 4)$, the $i$th path connecting
$x$ and $y$. Then, with respect to $N_u(e)$ and $N_v(e)$, there are
three cases to discuss.

\noindent {\bf Case $1$.} $x,y$ are in different sets. We claim that
$$
|n_u(e)-n_v(e)|=2|b_i-a_i|,
$$
where $a_i$ (resp. $b_i$) is the distance between $x$ (resp. $y$)
and the edge $e$.

To see this, assume that $x \in N_u(e), \ y \in N_v(e)$. Then we
have $a_i-b_i$ vertices more in $N_u(e)$ than in $N_v(e)$ on the
path $P_i$, but on each path $P_j$ $(j\neq i)$, we have $b_i-a_i$
vertices more in $N_u(e)$ than in $N_v(e)$. Hence
$|n_u(e)-n_v(e)|=|3(b_i-a_i)+(a_i-b_i)|=2|b_i-a_i|.$

\noindent{\bf Case $2$.} $x,y$ are in the same set. We claim that
$$
|n_u(e)-n_v(e)|=|V(G)|-g,
$$
where $g$ is the length of the shortest cycle of $G$ that contains
$e$.

To see this, assume that $x, y \in N_u(e)$. Thus all vertices from
the paths $P_j$ $(j\neq i)$ are in $N_u(e)$. Therefore,
$n_v(e)=\lfloor\frac{g}{2}\rfloor$, while
$n_u(e)=\lfloor\frac{g}{2}\rfloor+|V(G)|-g$. So
$|n_u(e)-n_v(e)|=|V(G)|-g.$

\noindent {\bf Case $3$.} One of $x,y$ is in $N_0(e)$. We claim that
$$
|n_u(e)-n_v(e)|\geq 2(a-1),
$$
with equality if and only if two paths of $P_i$ ($i=1,2,3,4$) have
length $a$, where $a$ is the length of a shortest path of the four
paths $P_i$ ($i=1,2,3,4$).

To see this, assume that $x \in N_u(e)$, $y \in N_0(e)$. Then the
shortest cycle $C$ of $G$ that contains $e$ is odd. Let $z_j \in P_j
(P_j \nsubseteq C)$ be the furthest vertex from $e$ such that $z_j \in
N_0(e)$. Then $|n_u(e)-n_v(e)|= \displaystyle\sum_{j}(d(x,z_j)-1) \geq
\displaystyle\sum_{j}(a+d(y,z_j)-1)\geq 2(a-1).$

From the above, we know that $|n_u(e)-n_v(e)|\geq 2$ in Case $2$. In
Case $3$, $|n_u(e)-n_v(e)| \leq 1$ if two paths of $P_i$ ($i=1,2,3,4$) have
length $1$, \ which is impossible since $G$ is simple. So,
$|n_u(e)-n_v(e)| \leq 1$ if and only if $x,y$ are in different sets and
$|b_i-a_i|=0$, that is, $e$ is in the middle position of an odd path
of $P_i$ ($i=1,2,3,4$).
\end{pf}
\begin{qed}
\end{qed}

\begin{lem}\label{lem3.2}
If $G$ is a  $\Theta_1$-graph of order $n\geq 12$. Then, we have
$$
Sz^*(G)<Sz^*(F_n)
$$
\end{lem}

\begin{pf}Without loss of generality, assume that $a\leq b\leq c \leq d$, then $b \geq 2.$ Now consider the six edges
which are incident with $x$ and $y$ but do not belong to  $P(a)$. Let
$e_1=xz$ be one of them, by Lemma \ref{lem3.1}, $\delta (e_1) \geq 2$. Similar thing is
true for the other five edges. Hence
$$
\displaystyle\sum_{e\in E(G)}\delta^2(e) \geq 6\times 2^2=24>18.
$$
Combining with equality $(1)$, this completes the proof.
\end{pf}
\begin{qed}
\end{qed}

\begin{lem}\label{lem3.3}
If $G$ is a  $\Theta_2$-graph of order $n\geq 12$. Then, we have
$$
Sz^*(G)<Sz^*(F_n)
$$
\end{lem}

\begin{pf}
Without loss of generality, let $d \geq b, e \geq c$.
In order to complete the proof, we consider the following
four cases.

\noindent {\bf Case $1$.} $d \geq b+2$.

Consider the two edges $xx_1, yy_1$ which belong to $P(d)$, then
$$
\delta (xx_1) =  \delta (yy_1)=\left\{
\begin{array}{ll}
a+c+e-2,&  b\leq a+c,\\
b+e-2, &   b\geq a+c.
\end{array}
\right.
$$

Therefrom, we get
$$\delta (xx_1) =  \delta (yy_1)\geq a+c+e-2.$$

Since $c+e\geq 3$, $a+c+e  \geq 4.$ If $a+c+e  \geq 6$,
then $\delta (xx_1) =  \delta (yy_1)\geq 4$, so
$\sum_{e\in E(G)}\delta^2(e) \geq 2\times 4^2>18.$

If $a=2,c=1,e=2,$ then $\delta (xx_1) =  \delta (yy_1)\geq 3$.
Now consider the edge $xx' \in P(e), \delta(xx') \geq 2$. So
$\sum_{e\in E(G)}\delta^2(e) \geq 2\times 3^2 + 2^2>18.$

If $a=1,c=1,e=3,$
since $n \geq 12, b+d-1 \geq 8$. Now consider the edge $xx' \in P(e), \delta(xx') \geq b+d-1 \geq 8$. So
$\sum_{e\in E(G)}\delta^2(e) \geq 8^2 >18.$

If $a=1,c=2,e=2,$ then $\delta (xx_1) =  \delta (yy_1)\geq 3$.
Now consider the edge $xx' \in P(e), \delta(xx') \geq 2$. So
$\sum_{e\in E(G)}\delta^2(e) \geq 2\times 3^2 + 2^2>18.$

If $a=1,c=1,e=2,$ if $b \geq 4>2=a+c,$ then $\delta (xx_1) =  \delta (yy_1)\geq 4$, so
$\sum_{e\in E(G)}\delta^2(e) \geq 2\times 4^2 >18.$ If $b=3$ or
$2, \delta (xx_1) =  \delta (yy_1)\geq 2, d\geq 7.$ Now consider the edge $zz' \in P(e), \delta(zz') \geq 4$. So
$\sum_{e\in E(G)}\delta^2(e) \geq 2\times 2^2 + 4^2>18.$
If $b=1$, then $d \geq 9$. Now consider the edge $xx' \in P(e), \delta(xx') \geq d \geq 9$. So
$\sum_{e\in E(G)}\delta^2(e) \geq 9^2>18.$

\noindent {\bf Case $2$.} $d = b+1, e=c+1$.

\noindent {\bf Subcase $2.1$.} $a+c-1 \geq b$.

Consider two edges $xx_1 \in P(c)$ and $xx_2 \in P(e)$,
$\delta (xx_1) \geq
d-1+e-2=b+e-2,$
$$
\delta (xx_2) = \left\{
\begin{array}{ll}
d+b-1,&  c\leq a+b,\\
d-1+c-1, &   c\geq a+b.
\end{array}
\right.
$$
Therefrom, we get $\delta (xx_2) \geq d+b-1=2b$. So,
$\delta^2(xx_1)+ \delta ^2(xx_2)= (b+e-2)^2+4b^2
=5b^2+ 2(e-1)b+(e-1)^2+3$.

If $b \geq 2$ or $e \geq 4$,
$\sum_{e\in E(G)}\delta^2(e) \geq \delta^2(xx_1)+ \delta ^2(xx_2)>18.$

If $b=1$, and $e \leq 3$, Now consider the edge $xx' \in P(d), \delta(xx') \geq 4$. So $\sum_{e\in E(G)}\delta^2(e) \geq 1^2 +2^2 +4^2>18.$

\noindent {\bf Subcase $2.2$.} $b \geq a+c+1$.

Consider the edge $xx_1 \in P(c)$, since $b \geq a+c+1,
 y\in N_{x_1}(xx_1)$. Let $u$ be the furthest vertex in $P(d)$ such that
 $u \in N_x(xx_1)$, $u'$ be the vertex incident with $u$ but not in  $N_x(xx_1)$. If the cycle $P(d) \cup P(c)\cup P(a)$ is even, then
 $d(u,x)=d(u',y)+a+c-1$, that is $d(u,x)-d(u',y)=a+c-1$.
 If the cycle $P(d) \cup P(c)\cup P(a)$ is odd, then
 $d(u,x)+1=d(u',y)+a+c-1$, that is $d(u,x)-(d(u',y)-1)=a+c-1$.
 So we have $\delta(xx_1)=e-2+a+c-1=a+2c-2$.

Then consider the edge $xx_2 \in P(e)$, since $b \geq a+c+1,
 y\in N_{x_2}(xx_2)$. Let $u_i (i=1,2)$ be the furthest vertex in $P(b)$
  and $P(d)$ such that
 $u_i \in N_x(xx_2)$, $u_i' (i=1,2)$ be the vertex incident with $u_i$ but not in  $N_x(xx_2)$. If the cycle $P(b) \cup P(c)\cup P(a)$ is even, then
 $d(u_1,x)=d(u_1',y)+a+c$, $d(u_2,x)+1=d(u_2',y)+a+c$.
 If the cycle $P(b) \cup P(c)\cup P(a)$ is odd, then
 $d(u_1,x)+1=d(u_1',y)+a+c$, $d(u_2,x)=d(u_2',y)+a+c$.
 So we have $\delta(xx_2)=d(u_1,x)+d(u_2,x)\geq 2a+2c-1$.

 From above, we have
 $$
 \sum_{e\in E(G)}\delta^2(e) \geq (a+2c-2)^2 +(2a+2c-1)^2 >18.
 $$
unless $a=c=1$. If $a=c=1$, now consider the edge $zz'$ belonging
to $P(e), \delta (zz') \geq 3$, so
$\sum_{e\in E(G)}\delta^2(e) \geq 1^2+3^2+3^2 >18$.

\noindent {\bf Subcase $2.3$.} $b=a+c$.

Consider the edge $xx_1 \in P(e)$, then
$\delta(xx_1)=d-1+b-1=2b-1$.

If $b \geq 3$,
then $\sum_{e\in E(G)}\delta^2(e) \geq 5^2 >18$.

If $b=2$, then $a=c=1,e=2,d=3$, which is impossible
since $n \geq 12$.

\noindent {\bf Case $3$.} $d=b+1, e=c$.

First, we know that $e=c \geq 2$.

\noindent {\bf Subcase $3.1$.} $a+c-1 \geq b$.

Consider the edges $xx_1 \in P(c)$ and $xx_2 \in P(e)$,
then
$$
\delta (xx_1)=\delta (xx_2)\geq d-1+e-1 =d+e-2.
$$
Since $d \geq 2$ and $e\geq2$, $d+e \geq 4$.

If $d+e \geq 6$, then $\sum_{e\in E(G)}\delta^2(e) \geq 2\times4^2 >18$.

If $4 \leq d+e \leq 5$, now consider the edge $xx' \in P(d)$.
If $d=3,e=2$, then $b=c=2, a \geq 5, \delta(xx') \geq 3$.
If $d=2,e=3$, then $b=1,c=3, a \geq 5, \delta(xx') \geq 5$.
If $d=2,e=2$, then $b=1,c=2, a \geq 7, \delta(xx') \geq 4$.
So $\sum_{e\in E(G)}\delta^2(e) >18$.

\noindent {\bf Subcase $3.2$.} $b>a+c-1$.

Consider the edge $xx_1 \in P(c)$, since $b >a+c-1$, then
$y \in N_{x_1}(xx_1)$. Let $u$ be the furthest vertex in $P(d)$
such that $z \in N_x(xx_1)$, $u'$ be the vertex incident with $u$
but not in $N_x(xx_1)$. If the cycle $P(d) \cup P(c)\cup P(a)$ is even, then
 $d(u,x)=d(u',y)+a+c-1$, $d(u,x)-d(u',y)=a+c-1$ .
 If the cycle $P(b) \cup P(c)\cup P(a)$ is odd, then
 $d(u,x)+1=d(u,y)+a+c-1$, $d(u,x)-(d(u',y)-1)=a+c-1$.
 So we have $\delta(xx_1)=(e-1) + (a+c-1)=a+2c-2$.

 Similarly
 $$
 \delta(xx_2)=a+2c-2.
 $$
where $xx_2$ is the edge belonging to $P(e)$.

Since $c \geq 2, a+2c \geq 5$.

If $a+2c \geq 6$, then $\sum_{e\in E(G)}\delta^2(e) \geq 2 \times 4^2>18$.

If $a+2c =5$, that is $a=1,c=e=2$, then $b \geq 4$. Now consider
$yy' \in P(d)$, then $\delta(yy') \geq 3$. So $\sum_{e\in E(G)}\delta^2(e) >18$.

\noindent {\bf Case $4$.} $d=b, e=c$.

\noindent {\bf Subcase $4.1$.} $b=d=c=e\geq 2$.

Consider the edge $xx_1 \in P(b)$, then $\delta(xx_1)=2(e-1)$.
Similarly for the other three edges incident with $x$.

If $e \geq 3$, then $\sum_{e\in E(G)}\delta^2(e) \geq 4 \times 4^2 >18$.

If $e =2$, since $n \geq 12, a\geq 6$. Now consider the edges $yy',zz'$
belonging to $P(a)$, $\delta (yy')=\delta(zz')\geq 2,$
so $\sum_{e\in E(G)}\delta^2(e) \geq 4 \times 2^2+2^2>18$.

\noindent {\bf Subcase $4.2$.} $b=d>c=e\geq 2$.

Consider the edge $xx_1 \in P(b)$, $\delta(xx_1)=d-1+e-1=d+e-2$.
For $xx_2 \in P(d)$, we also have $\delta(xx_2) =d+e-2.$

If $d+e\geq 6$, then $\sum_{e\in E(G)}\delta^2(e) \geq 2 \times 4^2 >18$.

If $d+e=5$, that is $d=3,e=2$, then $a \geq 4$. Now consider
$xx' \in P(c)$, then $\delta(xx') \geq 4$. So $\sum_{e\in E(G)}\delta^2(e) >18$.

Combining with equality $(1)$, this completes the proof.
\end{pf}
\begin{qed}
\end{qed}

\begin{lem}\label{lem3.4}
If $G$ is a  $\Theta_3$-graph of order $n\geq 12$. Then, we have
$$
Sz^*(G)<Sz^*(F_n)
$$
\end{lem}

\begin{pf}
Without loss of generality, let $f \geq d, e \geq c$.
In order to complete the proof, we consider the following
four cases.

\noindent {\bf Case $1$.} $e \geq c+2$.

Consider the edge $ww_1, yy_1 \in P(e)$,
$$
\delta(yy_1)=\delta(ww_1)=\left\{
\begin{array}{ll}
a+b+d+f-2, & c\leq a+b+d,\\
c+f-2, & c\geq a+b+d.
\end{array}
\right.
$$
Therefrom we get
$$
\delta(yy_1)=\delta(ww_1)\geq a+b+d+f-2.
$$
Since $d+f \geq 3, a+b+d+f \geq 5.$

If $a+b+d+f \geq 6$, then
$\sum_{e\in E(G)}\delta^2(e) \geq 2 \times 4^2 > 18$ .

If $a+b+d+f =5$, that is $a=b=d=1,f=2$. Now consider the edge $zz'\in P(f)$ then
$\delta(zz')\geq 2$, so
$\sum_{e\in E(G)}\delta^2(e) \geq 2 \times 3^2 +2^2 > 18$ .

\noindent {\bf Case $2$.} $e=c+1, f=d+1$.

\noindent {\bf Subcase $2.1$.} $a+c-1 \geq b+d$.

Consider the edge $yy_1 \in P(c)$, $yy_2 \in P(e)$, then
$
\delta(yy_1) =e-2+f-1=c+d-1,$
$$
\delta(yy_2)=\left\{
\begin{array}{ll}
b+d+f-1, & c\leq a+b+d,\\
c+f-2, & c\geq a+b+d.
\end{array}
\right.
$$
Therefrom, we get
$
\delta(yy_2) \geq b+d+f-1=b+2d.
$

If $d \geq 2$ or $b\geq 3$ or $c \geq 4$, then $\sum_{e\in E(G)}\delta^2(e) > 18$ .

If $d=1, b\leq 3, c\leq 3$, then consider the edge $xx' \in P(f)$, we have $\delta(xx') \geq 3$,
so $\sum_{e\in E(G)}\delta^2(e) \geq 1^2+ 3^2 +3^2 > 18$ .

\noindent {\bf Subcase $2.2$.} $a+c \leq b+d-1$.

It's similar to the Subcase 2.1.

\noindent {\bf Subcase $2.3$.} $a+c= b+d$.

Consider the edge $yy_1 \in P(e), xx_1 \in P(f)$, then
$\delta(yy_1)=b+d+f-2=b+2d-1,$
$\delta(xx_1)=a+c+e-2=a+2c-1.$
Since $n=a+b+c+d+e+f-2\geq 12$, then $(a+2c-1)+(b+2d-1) \geq 10$,
so $\sum_{e\in E(G)}\delta^2(e) \geq (a+2c-1)^2+ (b+2d-1)^2 > 18$ .

\noindent {\bf Case $3$.} $e=c+1, f=d$.

\noindent {\bf Subcase $3.1$.} $a+d-1 \geq b+c$.

Consider the edge $zz_1 \in P(d)$,
$
\delta(zz_1)\geq e-1+f-1=c+d-1.
$
Similarly
$\delta(zz_2)\geq c+d-1,$
where $zz_2$ is the edge belonging to $P(f)$.

Since $d \geq 2$, otherwise $G$ is not simple, then $c+d \geq 3.$

If $c+d\geq 5$, then $\sum_{e\in E(G)}\delta^2(e) \geq 2 \times 4^2> 18$.

If $c=1,d=3$, then $\delta(zz_1),\delta(zz_2) \geq 3$. Now consider the edge $yy' \in P(e),\delta(yy') \geq 3$,
so $\sum_{e\in E(G)}\delta^2(e) \geq 2 \times 3^2+3^2> 18$.

If $c=2,d=2$, then $\delta(zz_1),\delta(zz_2) \geq 3$. Now consider the edge $yy' \in P(e),\delta(yy') \geq 3$,
so $\sum_{e\in E(G)}\delta^2(e) \geq 2 \times 3^2+3^2> 18$.

If $c=1,d=2$, then $\delta(zz_1),\delta(zz_2) \geq 2$ and $e=f=2$. Now consider the edge $yy' \in P(e)$,
no matter $b \geq 2$ or $b=1$, we both have $\delta(yy') \geq 4$,
so $\sum_{e\in E(G)}\delta^2(e) \geq 2 \times 2^2+4^2> 18$.

\noindent {\bf Subcase $3.2$.} $a+d \leq b+c$.

Now consider the edge $ww_1 \in P(e)$, then
$$
\delta(ww_1)=\left\{
\begin{array}{ll}
a+d+f-2, & c\leq a+b+d,\\
c+f-2, & c\geq a+b+d.
\end{array}
\right.
$$
Therefrom, we get
$\delta(ww_1) =a+d-1+f-1=a+2d-2$.

Since $d \geq 2$, $a+2d \geq 5$.

If $a+2d \geq 7,$ then $\delta (ww_1) \geq 5$. So $\sum_{e\in E(G)}\delta^2(e) \geq 5^2> 18$.

If $a+2d =6,$ that is $a=2,d=2$, then $\delta (ww_1) \geq 4$. Now consider the edge $yy' \in P(e)$, $\delta(yy')\geq 2$. So $\sum_{e\in E(G)}\delta^2(e) \geq 4^2+2^2> 18$.

If $a+2d =5,$ that is $a=1,d=2$, then $\delta (ww_1) \geq 3$. Now consider the edge $yy' \in P(e)$, then we have  $\delta(yy')\geq \lceil{\frac{b+c+3}{2}}\rceil-1$. Since $n \geq 12,$ $b+2c\geq 8$. Then we have $b+c\geq 6$ unless $b=1,c=4$. When $b=1,c=4$, we can draw the graph exactly, we also have $\delta(yy')\geq 4$. So $\sum_{e\in E(G)}\delta^2(e) \geq 3^2+4^2> 18$.

\noindent {\bf Case $4$.} $d=f, e=c$.

We may assume that $a\leq b$.

\noindent {\bf Subcase $4.1$.} $c=e>d=f\geq 2$.

Consider the edge $ww_1 \in P(e)$, $\delta(ww_1)=f-1+c-1=c+f-2$.
For $ww_2 \in P(c)$, we also have $\delta(ww_2) =c+f-2.$

Since $c\geq 3$ and $f\geq 2$, $c+f\geq 5$.

If $c+f \geq 6$, then $\delta (ww_1) =\delta(ww_2)\geq 4$, so $\sum_{e\in E(G)}\delta^2(e) \geq 2 \times 4^2 >18$.

If $c+f =5,$ that is $c=3,f=2$, then $\delta (ww_1) =\delta(ww_2)\geq 3$. Now consider the edge $yy' \in P(e)$, then we have  $\delta(yy')\geq 1$. So $\sum_{e\in E(G)}\delta^2(e) \geq 2\times3^2+1^2> 18$.

\noindent {\bf Subcase $4.2$.} $c=e=d=f\geq 3$.

Consider the edge $ww_1 \in P(e), ww_2\in P(c)$, $\delta(ww_1)=\delta(ww_2) =f-1+c-1=2(c-1) \geq 4$. So $\sum_{e\in E(G)}\delta^2(e) \geq 2\times4^2> 18$.

\noindent {\bf Subcase $4.3$.} $c=e=d=f=2$.

If $b\geq a+4$, then we consider the edge $ww_1 \in P(e)$, $\delta(ww_1)=2$.
Similar for $ww_2\in P(c), xx_1\in P(d), xx_2 \in P(f)$. Then consider the edge $yy' \in P(b)$, $\delta(yy')\geq 2,$ so $\sum_{e\in E(G)}\delta^2(e) \geq 5\times2^2> 18$.

If $a\leq b\leq a+1$, then we consider the edge $ww_1 \in P(e)$, $\delta(ww_1)=2$.
Similar for $ww_2\in P(c), xx_1\in P(d), xx_2 \in P(f)$. Then consider the edge $yw_i,zx_i, (i=1,2)$, $\delta(yw_i)\geq 1,\delta(zx_i)\geq 1,$ so $\sum_{e\in E(G)}\delta^2(e) \geq 4\times2^2+4\times1^2> 18$.

If $b=a+3$, then we get $T_n$ with $n$ being odd.
If $b=a+2$, then we get $T_n$ with $n$ being even.

Combining with equality $(1)$, this completes the proof.
\end{pf}
\begin{qed}
\end{qed}

\begin{lem}\label{lem3.5}
If $G$ is a  $\Theta_4$-graph of order $n\geq 29$. Then, we have
$$
Sz^*(G)<Sz^*(F_n)
$$
\end{lem}

\begin{pf}
Without loss of generality, assume that $a=max\{a,b,c,d,e,f\}$.
Since $n \geq 29$, then $a \geq 6$. Now consider the edge
$ww_1 \in P(a)$. Then $z \in N_w(ww_1)$ or $z \in N_0(ww_1)$,
since $d(z,w) \leq d(z,w_1)$ by the choice of $a$.
And $z\in N_0(ww_1)$ if and only if $a=c\leq b+d$ and $e=1$.
We can obtain the similar result for $y$. Next, let $C$ be the shortest
cycle containing $ww_1$. Then $x \in N_w(ww_1)$, if $a > \frac{|C|+1}{2}$;
$x \in N_0(ww_1)$, if $a=\frac{|C|+1}{2}$; $x \in N_{w_1}(ww_1)$, if $a<\frac{|C|+1}{2}$.

\noindent {\bf Case $1$.} $a > \frac{|C|+1}{2}$.

Since $x \in N_w(ww_1)$, we can easily get $y,z \in N_w(ww_1)$.
So we have $\delta(ww_1)=n-|C|$. Similarly, $\delta(xx_1)=n-|C|$,
where $xx_1 \in P(a)$.

If $n-|C| \geq 4$, then $\sum_{e\in E(G)}\delta^2(e) \geq 2\times4^2> 18$.

If $n-|C| = 1$ and $C$ is composed of paths $P(a), P(f)$ and $P(b)$, then $V(G)-V(C)=\{z\}$,
and $e=c=d=1$. Since $P(a)\cup P(f) \cup P(b)$ is the shortest cycle, then
$f=b=1$ and $a\geq 26$, by $n\geq 29$. Now consider every edge $e$ in $P(a)$ except the middle one in P(a) when $a$ is odd, we have $\delta(e)=1$.
So $\sum_{e\in E(G)}\delta^2(e) \geq a-1> 18$.

If $n-|C| = 1$ and $C$ is composed of paths $P(a), P(f), P(d)$ and $P(c)$, which is impossible.

If $n-|C| = 2$ and $C$ is composed of paths $P(a), P(f)$ and $P(b)$, then $e+c+d \leq 4, f+b \leq 3$. Since $n\geq 29$, $a\geq 24$. Now consider the six edges $e_i(1\leq i\leq 6)$ in $P(a)$ such that the distance between $e_i$ and $x$
or $w$  no more than $2$, then we have $\delta(e_i)=2$.
So $\sum_{e\in E(G)}\delta^2(e) \geq 6\times2^2> 18$.

If $n-|C| = 2$ and $C$ is composed of paths $P(a), P(f), P(d)$ and $P(c)$, then one of the two vertices is in $P(b)$, another vertex is in $P(e)$. It is the case when  $C$ is composed of paths $P(a), P(f)$ and $P(b)$.

If $n-|C| = 3$ and $C$ is composed of paths $P(a), P(f)$ and $P(b)$, then $e+c+d \leq 5, f+b \leq 4$. Since $n\geq 29$, $a\geq 22$. Now consider the four edges $e_i(1\leq i\leq 4)$ in $P(a)$ such that the distance between $e_i$ and $x$
or $w$  no more than $1$, then we have $\delta(e_i)=3$.
So $\sum_{e\in E(G)}\delta^2(e) \geq 4\times3^2> 18$.

If $n-|C| = 3$ and $C$ is composed of paths $P(a), P(f), P(d)$ and $P(c)$, then either one of the two vertices in $P(b)$, another two vertices are in $P(e)$, or one of the two vertices in $P(e)$, another two vertices are in $P(b)$. It is the case when $C$ is composed of paths $P(a), P(f)$ and $P(b)$.

\noindent {\bf Case $2$.} $a = \frac{|C|+1}{2}$.

\noindent {\bf Subcase $2.1$.} $C$ is composed of paths $P(a), P(f), P(d)$ and $P(c)$.

In this case, $y,z \in N_w(ww_1)$ and $b>d+c$. Let $u$ be the furthest vertex in $P(e)$ such that $u\in N_w(ww_1)$, $u'$ be the vertex incident with $u$ but not in $N_w(ww_1)$. If the cycle $P(a)\cup P(c) \cup P(e)$ is even, then $d(x,u')+a-1=d(u,z)+c$, that is $d(u,z)=a-c-1+d(x,u')$. If the cycle $P(a)\cup P(c) \cup P(e)$ is odd, then $d(x,u')+a-1=d(u,z)+1+c$, that is $d(u,z)=a-c-2+d(x,u')$. Then $\delta(ww_1)=b-1+d(u,z)\geq a+b-c-3\geq a-1\geq 5$, since $b>d+c$. So $\sum_{e\in E(G)}\delta^2(e) \geq 5^2> 18$.

\noindent {\bf Subcase $2.2$.} $C$ is composed of paths $P(a), P(f)$ and $P(b)$.

In this case, $y \in N_w(ww_1)$ and $b\leq d+c$.

If $z \in N_0(ww_1)$, then $a=c \leq b+d$ and $e=1$. So
$\delta(ww_1)\geq c-1=a-1\geq 5$. Hence $\sum_{e\in E(G)}\delta^2(e) \geq 5^2> 18$.

If $z \in N_w(ww_1)$, similar to Subcase 2.1, we have
$$
d(u,z) \geq \left\{
\begin{array}{ll}
a-c-2,&  c\leq b+d,\\
a-(b+d)-2, &   c\geq b+d.
\end{array}
\right.
$$
Then $\delta(ww_1)=d-1+c+d(u,z)\geq a+d-3\geq a-2 \geq 4$.
Now consider the edge $xx_1 \in P(a)$. In this case, $w \in N_0(xx_1),
y \in N_x(xx_1)$. By the above analysis, if $z \in N_0(xx_1)$, then
$\delta(xx_1)\geq 5$. Hence $\sum_{e\in E(G)}\delta^2(e) \geq 5^2> 18$.
If $z \in N_x(xx_1)$, then
$\delta(xx_1)\geq 4$. Hence $\sum_{e\in E(G)}\delta^2(e) \geq 2\times 4^2> 18$.

\noindent {\bf Case $3$.} $a < \frac{|C|+1}{2}$.

\noindent {\bf Subcase $3.1$.} Both of $y$ and $z$ are in $N_0(ww_1)$.

In this case, $a=b=c, e=f=1$. Then $\delta(ww_1)=c-1=a-1\geq 5.$
Hence $\sum_{e\in E(G)}\delta^2(e) \geq 5^2> 18$.

\noindent {\bf Subcase $3.2$.} Both of $y$ and $z$ are in $N_w(ww_1)$.

In this case, we get
$$
\delta(ww_1) \geq \left\{
\begin{array}{ll}
a+d-2,&  d\geq |b-c|,\\
a+|b-c|-2, &   d\leq |b-c|.
\end{array}
\right.
$$
Then $\delta(ww_1)\geq a+d-2\geq a-1\geq 5.$
Hence $\sum_{e\in E(G)}\delta^2(e) \geq 5^2> 18$.

\noindent {\bf Subcase $3.3$.} One of $y$, $z$ is in $N_0(ww_1)$.

We may assume that $z\in N_0(ww_1)$, then $a=c\leq b+d, e=1$.

If $z \notin V(C)$, then $C=P(a)\cup P(f)\cup P(b)$. So
$\delta(ww_1)\geq c-1=a-1\geq 5$. Hence $\sum_{e\in E(G)}\delta^2(e) \geq 5^2> 18$.

If $z \in V(C)$, for $y \in N_w(ww_1)$, then $C=P(a)\cup P(e)\cup P(c)$.
Otherwise $C=P(a) \cup P(f) \cup P(d) \cup P(c)$, since $z \in N_0(ww_1)$,
then $y\in N_{w_1}(ww_1)$, a contradiction.
Let $u_1$ be the furthest vertex in $P(f)$ such that $u_1\in N_w(ww_1)$, $u_1'$ be the vertex incident with $u_1$ but not in $N_w(ww_1)$. If the cycle $P(a)\cup P(f) \cup P(b)$ is even, then $d(u_1,y)+b=d(u'_1,x)+a-1$, that is $d(u_1,y)-d(u_1',x)=a-b-1$. If the cycle $P(a)\cup P(f) \cup P(b)$ is odd, then $d(u_1,y)+b+1=d(u'_1,x)+a-1$, that is $d(u_1,y)-(d(u_1',x)-1)=a-b-1$.
Let $u_2$ be the furthest vertex in $P(d)$ such that $u_2\in N_w(ww_1)$, $u_2'$ be the vertex incident with $u_2$ but not in $N_w(ww_1)$. If the cycle $P(c)\cup P(e) \cup P(b)$ is even, then $d(u_2,y)+b=d(u'_2,z)+c=d(u'_2,z)+a$, that is $d(u_2,y)=a-b+d(u'_2,z)$. If the cycle $P(c)\cup P(e) \cup P(b)$ is odd, then $d(u_2,y)+b+1=d(u'_2,z)+a$, that is $d(u_2,y)=a-b-1+d(u'_2,z)$.
 Then $\delta(ww_1)=b+2(a-b-1)\geq 2a-b-2\geq a-2\geq 4$. Then consider the
 edge $xx_1$ in P(a), in this case, we have $w\in N_{x_1}(xx_1), z\in N_x(xx_1)$. If $y \in N_0(xx_1)$, by the above analysis, we have $\delta(xx_1) \geq 4$. So $\sum_{e\in E(G)}\delta^2(e) \geq a\times4^2> 18$. If $y\in N_x(xx_1)$, this is the Subcase 3.2.

Combining with equality $(1)$, this completes the proof.
\end{pf}
\begin{qed}
\end{qed}

From Lemma \ref{lem1}, \ref{lem2.1}, \ref{lem3.2}, \ref{lem3.3}, \ref{lem3.4} and \ref{lem3.5}, we have proved Theorem \ref{th1}.

\noindent {\bf Remark:} In fact, Theorem \ref{th1} can be improved
to $n \geq 23$, which needs more details of the proof. But
$n$ can not be decrease, because the revised Szeged index of the graph $\Theta_4$ with $b=c=d=e=f=1$ is less than $F_n$.

\end{document}